# Research into the Group of Two-dimensional Magic Cube and its Cayley Graph Diameter

Xinchang High School, Zhejiang, China

Author: Zihan Jia, Qixuan Zhang, Yuming Xiang

**Instructor: Yudong Wu**

**Zhejiang, China**

# Abstract

Based on the rules of magic cubes, a game of two-dimensional magic cube was deliberately designed. This essay will explore its properties with the assistance of group theory and computer programming. It will first elaborate the rules of two-dimensional magic cube and then use group theory to comprehensively study its properties like the order of the permutation group of the cube, the diameter of Cayley Graph, etc. At the last part of the paper, a much more general result will be raised to satisfy situations like the irregular magic cubes.

**Keywords:** two-dimensional magic cube; permutation group; the diameter of Cayley Graph

# Contents

# 1 Introduction

## 1.1 Foreword

The magic cube (Figure 1) has been a popular toy for a long time since it was invented by Hungary architect, Rubik in 1974. After that, various kinds of magic cubes were designed (Figure 2). However, a majority of the changes among them were only focusing on their orders and structures. Some fans attempted to compute a four-dimensional magic cubes game but failed to promote it to the public, due to its complicated controls, unclear models and inability to be replicated in real life.

Inspired by this, a new game of two-dimensional magic cube was raised by us based on the rules of standard magic cubes, which will be analyzed in the coming sections. Besides, it should be noted that the two-dimensional magic cube is not the projection of three-dimensional magic cubes but a game shared with similar rules.

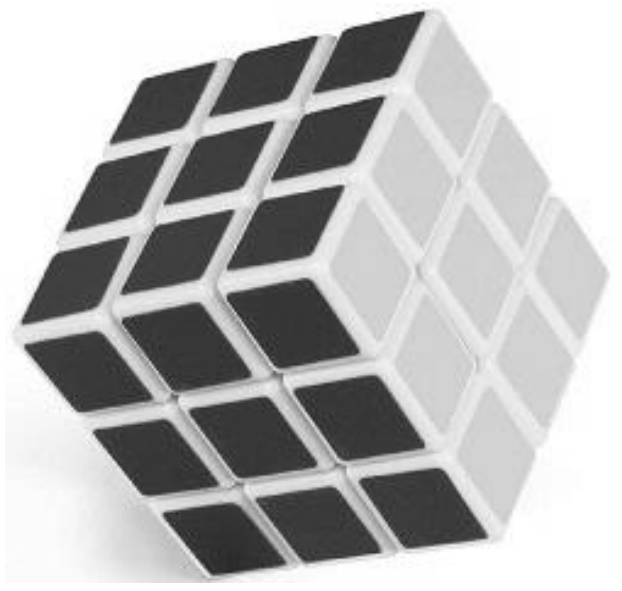

**Figure 1**

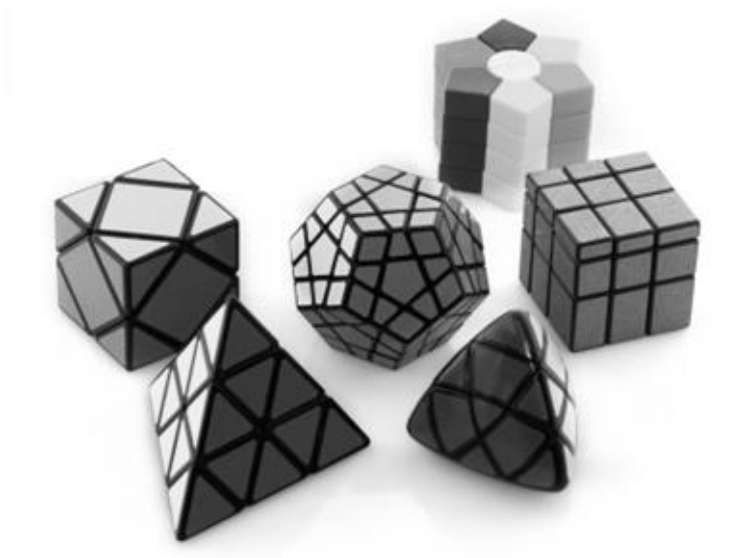

**Figure 2**

## 1.2 Introduction to the game rules

This section will elaborate the rules of two-dimensional magic cubes.

### 1.2.1 The structure of two-dimensional magic cubes

An $n$-order two-dimensional magic cube consists of $n^2$ congruent squares, with numbers similar to 0, 1, 2, …, $n^2$-1 on them. Each of the squares can be put at any position in the board, forming a larger square of $n \times n$ size (Figure 3).

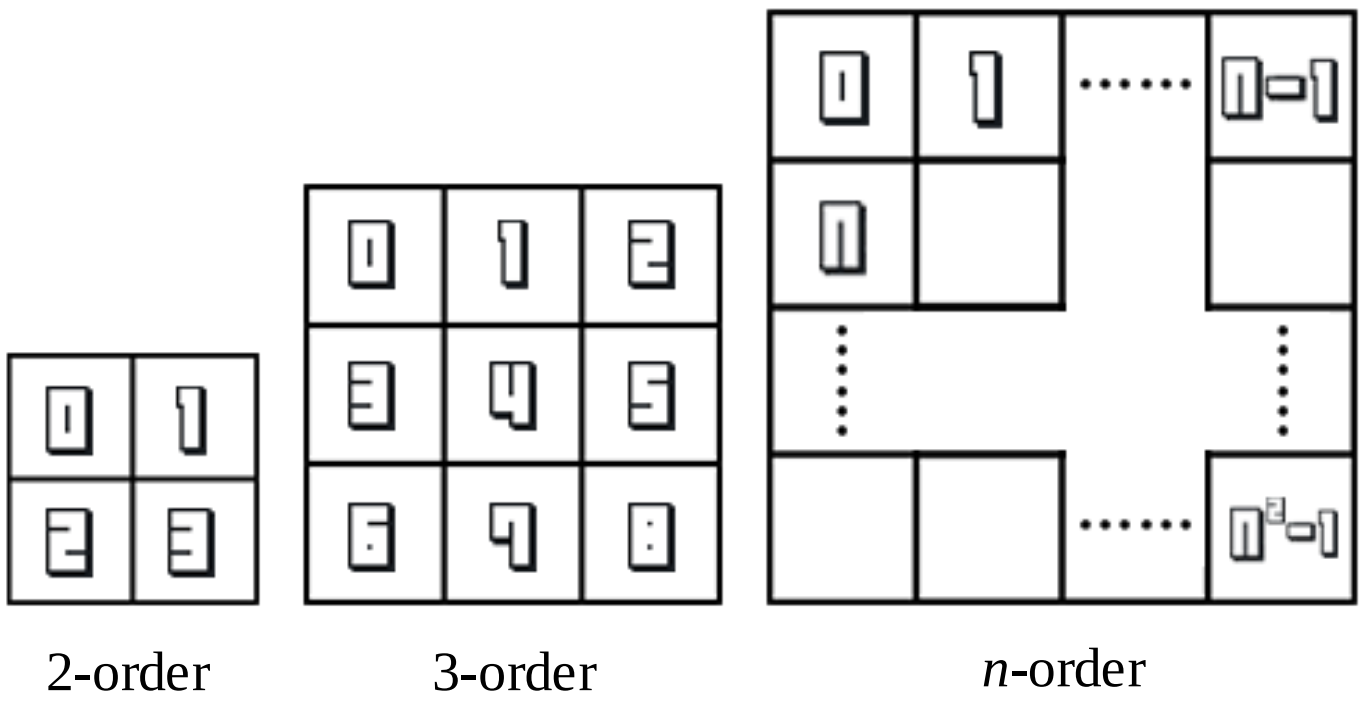

**Figure 3**

### 1.2.2 Manipulations

***Definition* 1.1** **Move** stands for the unit movement upwards or leftwards of any column or row of an *n*-order two-dimensional magic cube. The magic cube block beyond the square is then moved to the end of the column or row. **Manipulation** stands for moving certain column or row for any time except the multiples of *n*.

Figure 4 is a manipulation of a 3-order two-dimensional magic cube.

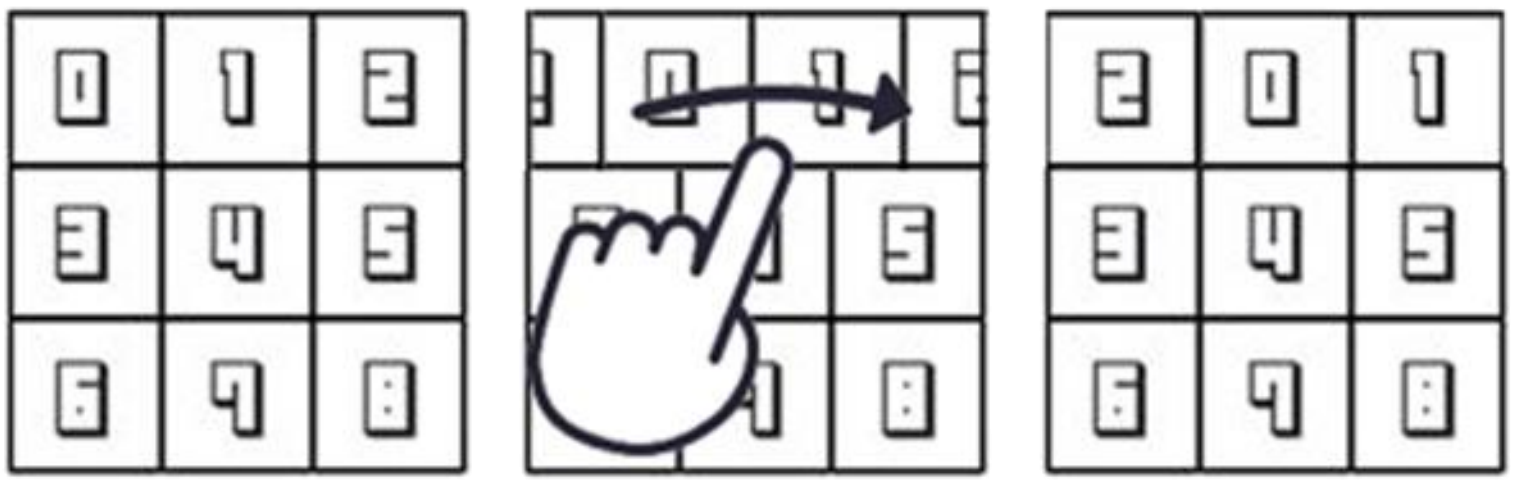

Before After

**Figure 4**

The results show that no matter how the *n*-order two-dimensional magic cube is manipulated, its size stays as an $n \times n$ square, filled with the same numbers: 0, 1, 2, …,$n^2$-1.

### 1.2.3 Initial state, solve and solution step number

Next, we will give the definitions of initial state, solve and solution step number.

***Definition* 1.2** A two-dimensional magic cube is in the **initial state**, if and only if the numbers of the cube is increasing from left to right, top to bottom.

***Definition* 1.3** A two-dimensional magic cube is **solved**, when it is transformed into the initial state from a not initial state after several manipulations. The number of the manipulations is called **solution step number**.

## 1.3 Proposed problems

People have put forward and solved many interesting mathematical problems over Rubik’s Cube as follows:

**Question 1** How many ways are there to scramble a Rubik’s Cube?

**Answer 1** There are 43, 252, 003, 274, 489, 856, 000 ways to scramble a Rubik’s Cube.

**Question 2** How many manipulations are needed at least to solve any Rubik’s Cube?

**Answer 2** We can solve any 3-order magic cube in 20 manipulations.

The answer to Question 2 is called the God’s number, also the Cayley Graph diameter of the Rubik’s Cube cube and it has troubled mathematicians for over three decades, until July 2010. California scientists expropriated more powerful resources, a supercomputer at Google headquarters in San Francisco. With the simplicity of the program and the improvement of the device, the question was solved by the super configured computer. The scientists used this computer, validated every single state by using the enumeration method, provided that any Rubik’s Cube cube could be solved in 20 manipulations, so 20 is the official God’s number.

What needs to be pointed out is that the Group Theory will be a powerful mathematical tool to solve problems of the magic cubes. The introduction of the Group Theory will be given in Section 2.1.

As regard to the two-dimensional magic cube, there are some questions as follows:

**Question 3** How many ways are there to scramble a 3-order two-dimensional magic cube?

**Question 4** How many ways are there to scramble an $n$-order two-dimensional magic cube?

**Question 5** How many manipulations are needed, at least, to solve any 3-order magic cube?

**Question 6** How many manipulations are needed, at least, to solve any $n$-order magic cube?

The answer to the Question 3, Question 4 and Question 5 will be given in the future sections. Due to limited abilities, we cannot solve the exact answer to the Question 6, but the lower bound is given by using estimations.

# 2 Model establishment

## 2.1 Fundamental theorem of group theory

***Definition* 2.1.1** A **group** is a set, $G$, together with an operation "·" (called the group law of $G$) that combines any two elements $a$ and $b$ to form another element denoted $a \cdot b$ or $ab$. To qualify as a group, the set and operation must satisfy four requirements known as the group axioms:
(1)Closure For all $a$, $b$ in $G$, the result of the operation, $a \cdot b$, is also in $G$.
(2)Associativity For all $a$, $b$ and $c$ in $G$, $(a \cdot b) \cdot c = a \cdot (b \cdot c)$.
(3)Identity element For each $a \in G$, $\exists e \in G$ satisfies $e \cdot a = a \cdot e = a$.
(4)Inverse element For each $a$ in $G$, $\exists a^{-1} \in G$, $a^{-1} \cdot a = a \cdot a^{-1} = e$.

***Definition* 2.1.2** A **finite group** is a mathematical group with a finite number of elements. The number of elements, its order, is commonly denoted as $|G|$.

***Definition* 2.1.3** Let $\Omega$ be a set including $n$ elements:

$$\Omega = \{\alpha_1, \alpha_2, \ldots, \alpha_n\}.$$

A **permutation** or an ***n*-permutation** is a bijection from $\Omega$ to itself. $\alpha_i^\sigma$ indicates the image of $\alpha_i$ under permutation $\sigma$. Define the multiplication of $\sigma$ and $\tau$ as continuous applying of them. That is to say:

$$i^{\sigma\tau} = (i^\sigma)^\tau, i = 1, 2, \cdots, n.$$

***Definition* 2.1.4** A **permutation group** is a group $G$ whose elements are permutations of a given set $M$ and whose group operation is the composition of permutations in $G$.

***Theorem* 2.1.5** The **symmetry group** of an object is the group of all permutations under which the object is invariant with composition as the group operation, commonly denoted as $S_n$. $|S_n|=n!$.

***Proof.*** [1]P11

***Definition* 2.1.6** A permutation is called a **cyclic permutation** if and only if it has a single nontrivial cycle. The cycle can be written by using the compact cycle notation $\sigma = (\alpha_{i_1}\ \alpha_{i_2}\ \dots\ \alpha_{i_m})$ (there are no commas between elements in this notation, to avoid confusion with an $m$-tuple). The length of a cycle is the number of elements. A cycle of length $m$ is also called an ***m*-cycle**. A cycle with only two elements is called a **transposition**.

***Theorem* 2.1.7** Arbitrary permutation $\sigma$ can be expressed as a product of a number of transpositions. The parity of the number of the transpositions is uniquely determined by $\sigma$. The parity of the $n$-arrangement $1^\sigma, 2^\sigma, \cdots, n^\sigma$ is the same as the parity of the number of transpositions.

***Proof.*** [1]P16

***Definition* 2.1.8** An **odd permutation** is a permutation obtainable from an odd number of transpositions. An **even permutation** is a permutation obtainable from an even number of transpositions.

***Theorem* 2.1.9** The product of any two even permutations is an even permutation.

***Proof.*** According to *Definition* 2.1.8, we can easily draw this conclusion.

***Theorem* 2.1.10** The **alternating group** of an object is the group of all even permutations under which the object is invariant with composition as the group operation, commonly denoted $A_n$. $|A_n|=n!/2$.

***Proof.*** [1]P17

***Theorem* 2.1.11** (**Lagrange's Group Theorem**) For any finite group $G$, the order (number of elements) of every subgroup $H$ of $G$ divides the order of $G$.

***Proof.*** [1]P29

***Theorem* 2.1.12** $A_n (n \geq 3)$ is generated by all cycles of length 3.

***Proof.*** [1]P77

***Definition* 2.1.13** In mathematics, given two groups, $(G, *)$ and $(H, \cdot)$, a **group homomorphism** from $(G, *)$ to $(H, \cdot)$ is a function $f$: $G \to H$ such that for all $x$ and $y$ in $G$ it holds that

$$f(x * y) = f(x) \cdot f(y),$$

An **isomorphism** is a bijective homomorphism. In this case, the groups $G$ and $H$ are called isomorphic. They differ only in the notation of their elements and are identical for all practical purposes.

## 2.2 The group of two-dimensional magic cube

We mark each position with number 1, 2, …, $n^2$, from top to bottom, left to right. It can also be described by a coordinate like $(a, b)$, which represents row $a$, column $b$.

Let $F_n=\{(x, y)\,|\,x, y\in\{1, 2, \cdots, n\}\}$, $B_n=\{1, 2, \cdots, n^2\}$, structure a map $\varphi$ from $F_n$ to $B_n$:

$$\varphi(x, y) = n(x-1) + y,\ (x, y)\in F_n,$$

Obviously, $\varphi$ is a one to one correspondence from the coordinate representation of position to number representation of position.

For each two-dimensional magic cube, we also have two kinds of methods to describe its state:

1. **Number-notation:** $(x_1, x_2, \cdots, x_{n^2})$, where $x_i$ denotes the number on position $i$.
2. **Position-notation:** $(y_1, y_2, \cdots, y_{n^2})$, where $y_i$ denotes the ordinal number of the position with the number $i$.

Let $C_n$ be the set of all permutations of all elements in $B_n$. Both number-notation and position-notation are permutations of all elements in $B_n$. Each two-dimensional magic cube corresponds to one number-notation and one position-notation. This paper will use position-notation in research of the group of the two-dimensional magic cube and number-notation in research of its Cayley Graph diameter.

Let $V_i^j$ represent a $j$-time move which is applied to column $i$ and $H_i^j$ represent a $j$-time move which is applied to row $i$. A positive number means moving upwards or leftwards, while a negative number means moving downwards or rightwards. Each manipulation can be described by $V_i^j$ or $H_i^j$.

Let $X = \{V_i^j\,|\,i = 1, 2, \cdots, n, j\in Z\} \cup \{H_i^j\,|\,i = 1, 2, \cdots, n, j\in Z\}$.

***Definition* 2.2.1 Transformation** is an $n^2$-permutation which is applied to $B_n$. Applying a transformation on a two-dimensional magic cube means applying the permutation to its position-notation.

$V_i^j, H_i^j$ are all transformations. For example, the number on position (2, 1) is moved to position (1, 1) after manipulation $V_1^1$. In other words, transformation $V_1^1$ maps position (2, 1) to position (1, 1).

Each element in $X$ can be denoted as follows:

$$V_i^j = (\varphi(n, i)\ \varphi(n-1, i)\ \varphi(n-2, i) \cdots \varphi(1, i))^j,$$

or

$$H_i^j=(\varphi(i, n)\ \varphi(i, n-1)\ \varphi(i, n-2) \cdots \varphi(i, 1))^j,$$

The composite operation between transformations follows the rules of permutation (from left to right). The composite operation combines any two transformations to form another transformation. We do not distinguish transformations with the same effect. For example, $V_1^1$ and $V_1^{n+1}$ are the same transformations and $V_1^1 V_2^1$ and $V_2^1 V_1^1$ are the same transformations. Let $G_n$ represent a set of all the transformations of $n$-order two-dimensional magic cubes, then any element in $G_n$ is a permutation applied to $B_n$.

***Theorem* 2.2.2** $G_n$ is a permutation group under composite operation, which is called the group of a two-dimensional magic cube.

***Proof.*** Every element in $G_n$ can be represented by a composition of some elements in $X$. Obviously for $\forall a, b \in G_n$, $\exists c = ab \in G_n$.

(1) According to the rules of the two-dimensional magic cube, $\forall M_1, M_2, M_3 \in G_n$ satisfy:

$$M_1M_2M_3 = (M_1M_2)M_3 = M_1(M_2M_3).$$

(2) Identity permutation $e$ maps $M_0$ to $M_0$, where $e$ is the identity element.

(3) For $\forall M \in G_n$, $M$ can be denoted as follows:

$$M = M_1M_2 \cdots M_m, M_i \in X, i \in \{1, 2, \cdots, m\},$$

then

$$M^{-1} = (M_m)^{-1} \cdots (M_2)^{-1} (M_1)^{-1}.$$

And for $\forall V_i^j, H_i^j \in X$,

$$\left(V_i^j\right)^{-1} = V_i^{-j} \in X, \left(H_i^j\right)^{-1} = H_i^{-j} \in X,$$

which means $\forall M \in G_n$, $\exists M' = M^{-1} \in G_n$, $M^{-1}$ is the inverse element of $M$.
According to *Definition* 2.1.1, we can draw the conclusion.

# 3 The order of $G_n$

In this chapter, we will prove $|G_n| = \begin{cases} n^2!/2, n \equiv 1 \pmod 2, \\ n^2!, n \equiv 0 \pmod 2. \end{cases}$

***Lemma* 3.1** $\forall n \geq 3, n \in Z, \exists M_0 \in G_n, M_0 = (1\ 2\ 3)$.

***Proof.*** It is easy to verify that $M_0 = V_1^{-1} H_1^1 V_1^1 H_1^1 V_1^{-1} H_1^{-2} V_1^1 \in G_n$ and $M_0 = (1\ 2\ 3)$.

***Lemma* 3.2** $\forall n \geq 3, n \in Z, \{a, b, c\} \subseteq B_n, \exists M_0 \in G_n, M_0 = (a\ b\ c)$.

***Proof.*** (1) If $a = \varphi(1, y)$, then let $M_1 = H_1^{y-1}$, else $a = \varphi(x, y)(x \neq 1)$, then let $M_1 = H_1^{y-1} V_1^{x-1}$, which makes $a^{M_1} = \varphi(1, 1) = 1$.

(2) If $b^{M_1} = \varphi(1, y)$, then let $M_2 = V_1^{-1} H_1^{y-2} V_1^1$, else $b^{M_1} = \varphi(x, y)(x \neq 1)$, then let $M_2 = H_1^{y-2} V_1^{x-1}$, which makes $b^{M_1M_2} = \varphi(1, 2) = 2$, and $a^{M_1M_2} = 1^{M_1} = 1$.

(3) If $c^{M_1M_2} = \varphi(1, y)$, then let $M_3 = V_1^{-1} H_1^{y-3} V_1^1$, else $c^{M_1M_2} = \varphi(x, y)(x \neq 1)$, then let $M_3 = H_1^{y-3} V_1^{x-1}$, which makes $c^{M_1M_2M_3} = \varphi(1, 3) = 3$, and $a^{M_1M_2M_3} = 1^{M_3} = 1$, $b^{M_1M_2M_3} = 2^{M_3} = 2$.

Let $M = M_1M_2M_3$, then $a^M = 1, b^M = 2, c^M = 3$. From *Lemma* 3.1, we can know that there is $M' \in G_n$, $M' = (1\ 2\ 3)$. Considering $M_0 = MM'M^{-1}$, after the transformation $M_0$ to $a, b, c$, we can get $a^{M_0} = a^{MM'M^{-1}} = 1^{M'M^{-1}} = 2^{M^{-1}} = b$. Similarly, $b^{MM'M^{-1}} = c$, $c^{MM'M^{-1}} = a$. But for any other number $p$, it is obvious that $p^M \neq 1, 2, 3$, $p^{MM'} = p^M$, $p^{M_0} = p^{MM'M^{-1}} = p^{MM^{-1}} = p$. Thus, $M_0 = (a\ b\ c)$. Proved.

***Theorem* 3.3** $G_n=\begin{cases}A_{n^2}, n\equiv 1 \pmod 2,\\ S_{n^2}, n\equiv 0 \pmod 2.\end{cases}$

***Proof.*** As $G_n$ is a permutation group applied to $B_n=\{1, 2, \cdots, n^2\}$, therefore, $G_n\leq S_{n^2}$.

1. When $n = 2$, then $G_2 \leq S_4$, $H_1^1 = (1,2)$, $V_1^1=(1,3)$, $H_1^1V_2^1H_1^1=(1,4)$, $H_1^1V_1^1H_1^1=(2,3)$, $V_2^1=(2,4)$, $H_2^1=(3,4)$.

   According to *Theorem* 2.1.7, such permutations can generate all the 4-cycles in 2-order two-dimensional magic cube, which means $G_2=S_4=S_{2^2}$.

2. When $n\geq 3$, by applying *Lemma* 3.2, it is easy to find all the 3-cycles applied to $B_n$ in $G_n$. Meanwhile, considering *Theorem* 2.1.12, all the 3-cycles can generate $A_{n^2}$. Therefore, $G_n\geq A_{n^2}$.

   (1) If $n\equiv 1 \pmod 2$, then $G_n$ is generated by $n$-element even permutations $H_1^1$, $H_2^1, \cdots$, $H_n^1$, $V_1^1$, $V_2^1, \cdots$, $V_n^1$. According to *Theorem* 2.1.9, there is no odd permutation in $G_n$. Hence, $G_n\leq A_{n^2}$, which means $G_n=A_{n^2}$.

   (2) If $n\equiv 0 \pmod 2$, as odd permutation $H_1^1\in G_n$, then $G_n\neq A_{n^2}$, which means $G_n>A_{n^2}$, $|G_n|>|A_{n^2}|=n^2!/2$. Also, $G_n\leq S_{n^2}$, according to *Theorem* 2.1.11, then $|G_n|\,|\,|S_{n^2}|=n^2!$. Thus, $|G_n|=n^2!$, $G_n=S_{n^2}$.

***Theorem* 3.4** $|G_n|=\begin{cases}n^2!/2, n\equiv 1 \pmod 2,\\ n^2!, n\equiv 0 \pmod 2.\end{cases}$

***Proof.*** By applying *Theorem* 2.1.5, *Theorem* 2.1.10 and *Theorem* 3.3, we can draw the conclusion.

# 4 The Cayley Graph diameters

The diameter of the Cayley Graph of a two-dimensional magic cube has been researched in this section. Some questions in this section has been referenced to [4].

## 4.1 Relevant definitions

***Definition* 4.1.1** A **graph** is an ordered pair $G = (V, E)$ comprising a set $V$ of vertices, together with a set $E$ of edges which are 2-element subsets of $V$.

***Definition* 4.1.2** The **degree** of the vertex $v$ is the number of edges that connect to it, commonly denoted as deg($v$).

***Definition* 4.1.3** A **path graph** of order $n \geq 2$ is a graph in which the vertices can be listed in an order $v_1, v_2, \cdots, v_n$ such that the edges are the $\{v_i, v_{i+1}\}$ where $i = 1, 2, \ldots, n-1$.

***Definition* 4.1.4** In an undirected graph, an unordered pair of vertices $\{x, y\}$ is called **connected** if a path leads from $x$ to $y$. A **connected graph** is an undirected graph in which every unordered pair of vertices in the graph is connected.

***Definition* 4.1.5** The **distance** between two vertices is the length of the shortest path between those vertices. The **diameter** $d$ of a graph is the maximum eccentricity of any vertex in the graph. That is, $d$ is the greatest distance between any pair of vertices or, alternatively, $\operatorname{diam}((V, E))=\max\{d(v, w) \mid v, w \in V\}$.

***Definition* 4.1.6** Let $G = \langle g_1, g_2, \cdots g_n \rangle$, then $G$ is a permutation group generated by set $X$, where $X=\{g_1, g_2, \cdots, g_n\}$. The Cayley Graph of $G$ on $X$ is a graph $(V, E)$, where $V$ consists of all the elements in $G$. The edges of Cayley Graph satisfy: if $x, y \in V = G$, then $x$ and $y$ are connected by an edge if and only if $y = g_i x$ or $x = g_i y$, $i = 1, 2, \cdots, n$.

***Theorem* 4.1.7** Let $\Gamma_G = (V, E)$ denote the Cayley Graph of permutation group $G = \langle g_1, g_2, \cdots, g_n \rangle$. For $\forall v \in V$, $\deg(v)=\left|\{g_1, g_2, \cdots, g_n\} \cup \{g_1^{-1}, g_2^{-1}, \cdots, g_n^{-1}\}\right|$.

***Proof.*** [4]P13.

According to *Theorem* 4.1.7 , for any $v$ in two-dimensional magic cube group $G_n = \langle X \rangle$ ,deg($v$) = $|X| = 2n(n-1)$.

***Theorem* 4.1.8** Let $\Gamma_G = (V, E)$ denote the Cayley Graph of a two-dimensional magic cube group $G_n$ . For $\forall v,w \in V$, $\exists u \in V$ which satisfies $d(v, w) = d(e, u)$, where $e$ represents the identity element.

***Proof.*** Without loss of generality, we can take $d(v, w) = x$, then $v g_1 g_2 \cdots g_x = w$, $g_i \in X$, $i = 1, 2, \cdots, x$. Let $u = v^{-1}w$ and it satisfies $d(v, w) = d(e, u)$.

***Corollary* 4.1.9** $\operatorname{diam}(\Gamma_G) = \max\{d(e, v), v \in V\}$.

***Proof.*** With *Definition* 4.1.5, together with *Theorem* 4.1.8 and *Theorem* 4.1.7, we get the *Corollary* 4.1.9.

In view of *Corollary* 4.1.9, the diameter of the Cayley Graph of a two-dimensional magic cube group is equal to the maximum distance between the identity element and any other element. Denote the set of the element from which the distance to the identity element is $i$ in the two-dimensional magic cube group $G_n$ by $G_n(i)$. Then

$$G_n(i) = \{v \in G_n \mid d(e, v) = i\}.$$

According to the set theory, the following *Theorem* 4.1.10-4.1.12 hold true.

***Theorem* 4.1.10** $G_n(i) \cap G_n(j) = \emptyset$, $i \neq j$.

***Theorem* 4.1.11** $G_n = \bigcup_{i=0}^{\infty} G_n(i)$ .

***Theorem* 4.1.12** $|G_n| = \sum_{i=0}^{\infty} |G_n(i)|$ .

***Theorem* 4.1.13** $\forall M \in X$, $v \in G_n(i)$, $w = vM$, then $w \in G_n(i-1) \cup G_n(i) \cup G_n(i+1)$.

***Proof.*** For $w = vM$, $v = wM^{-1}$,then$\left|d(e, w) - d(e, v)\right| \leq 1$, that is $d(e, v) - 1 \leq d(e, w) \leq d(e, v) + 1$. So $w \in G_n(i-1) \cup G_n(i) \cup G_n(i+1)$. The proof of *Theorem* 4.1.13 is completed.

So we can get the Cayley Graph diameter by finding an $i$ satisfies $|G_n(i)| \neq 0$ and $|G_n(i+1)| = 0$. Then $i$ is the diameter of $G_n$'s Cayley Graph .

## 4.2 God's algorithm

### 4.2.1 Basic logic of the algorithm

Structure a mapping $f : G_n \rightarrow C_n$, $f(v) = v(c)$, where $c$ denotes the number-notation of the two-dimensional magic cube in the initial state and $v(c)$ denotes the number-notation of the two-dimensional magic cube generated by applying $v$ to the initial magic cube. Obviously, $f$ is a bijection and $f(e) = c$.

Let $P_n = \{v(c)|v \in G_n\}$, then $P_n$ contains all the possible number-notations of $n$-order two-dimensional magic cube. Let "$\cdot$" be a binary operations on $P_n$ which satisfies the following relationship:

$$p \cdot q = (yw)(c), p, q \in P_n,$$

where $p = v(c)$, $q = w(c)$. Obviously $P_n$ is a group.

***Lemma* 4.2.1** $P_n \cong G_n$.

***Proof.*** Since that

$$f(v)f(w) = v(c)w(c) = (vw)(c) = f(vw), \forall v, w \in G_n,$$

and $f$ is a bijection, given *Theorem* 2.1.13, we can draw the conclusion that $P_n \cong G_n$.

Let $Y$ denote the generating set of $P_n$, then $Y = \{f(v)|v \in X\}$ also the image of $X$ under map $f$. It is easy to prove the following lemma:

***Lemma* 4.2.2** The Cayley Graph of $P_n$ is isomorphic to the Cayley Graph of $G_n$.

Define $P_n(i)$ by

$$P_n(i)=\{f(v)|v \in G_n(i)\},$$

Then $P_n(i)$ is the image of $G_n(i)$ under map $f$. Thus, $|G_n(i)| = |P_n(i)|$.

By the isomorphic relation we can obtain that $P_n(i) = \{p \in P_n | d(c, p) = i\}$.

We can also get the following lemmas:

***Lemma* 4.2.3** $P_n(i) \cap P_n(j) = \emptyset, i \neq j.$

***Lemma* 4.2.4** $P_n = \cup_{i=0}^{\infty} P_n(i)$.

***Lemma* 4.2.5** $|P_n| = \cup_{i=0}^{\infty} |P_n(i)|$.

***Lemma* 4.2.6** $\forall M \in X, v \in P_n(i), w = vM \in P_n(i\text{-}1) \cup P_n(i) \cup P_n(i+1)$.

***Lemma* 4.2.7** Let $G_n = \langle X \rangle, P_n = \langle Y \rangle$, then the diameter of $G_n$'s Cayley Graph is equal to the diameter of $P_n$'s Cayley Graph.

To get the Cayley Graph diameter of the two-dimensional magic cube, we enumerate all possibilities. The algorithm is shown as follows:

1. $i = 0$,
2. $P_n(0) = \{c\}$,
3. Take out an element $p$ from $P_n(i)$,
4. Multiply $p$ and every elements in $Y$ and we can get $n(n\text{-}1)$ elements $q_1, q_2, \cdots, q_{n(n-1)}$.

5. For each element $q$ in $q_1$, $q_2$, $\cdots$,$q_{n(n-1)}$, if $q \notin P_n(i\text{-}1) \cup P_n(i) \cup P_n(i+1)$, then put it into $P_n(i+1)$,

6. Repeat steps 3, 4 and 5, till all the elements in $P_n(i)$ are traversed,

7. If $|P_n(i)| \neq 0$, $i = i+1$, then go to step 3, else go to step 8,

8. The Cayley Graph diameter of $P_n$ is $i$.

9. Given *Theorem* 4.2.7, the Cayley Graph diameter of $G_n$ is $i$.

### 4.2.2 Program result of *n*=3

According to *Theorem* 3.4, $|G_3|=3^2!/2=181440$. We have written a program (codes in *Appendix* 1) to solve the Cayley Graph diameter of $G_3$. Using *Python 2.7*, we can get the output as follows:

depth: 1
the number of all states: 13
depth: 2
the number of all states: 109
depth: 3
the number of all states: 845
depth: 4
the number of all states: 6053
depth: 5
the number of all states: 34727
depth: 6
the number of all states: 124224
depth: 7
the number of all states: 178965
depth: 8
the number of all states: 181440
depth: 9
the number of all states: 181440
The God's number of two-dimensional magic cube 3*3 is 8.

Through the program we can also get the values of $\sum_{k=0}^{i} G_3(k)$, which are listed as follows:

**Table 1 Data of the Cayley Graph diameter of 3-order two-dimensional magic cube**

| $i$ | $\sum_{k=0}^{i} G_3(k)$ | $G_3(i)$ |
|---|---|---|
| 0 | 1 | 1 |
| 1 | 13 | 12 |
| 2 | 109 | 96 |
| 3 | 845 | 736 |
| 4 | 6053 | 5208 |

**Continued table**

| $i$ | $\sum_{k=0}^{i} G_3(k)$ | $G_3(i)$ |
|---|---|---|
| 5 | 34727 | 28674 |
| 6 | 124224 | 89497 |
| 7 | 178965 | 54741 |
| 8 | 181440 | 2475 |

## 4.3 Estimates of the lower bound

In the practical application of the God's algorithm, we found that with the increasing of the order, the number of total states of two-dimensional magic cube is growing so rapidly that the Cayley Graph diameters can hardly be solved by using a personal computer. Therefore, we use the mathematical method to estimate the lower bounds of the Cayley Graph diameters of two-dimensional magic cube of high orders.

### 4.3.1 Theoretical analysis

We estimate the lower bound of the Cayley Graph diameter of a two-dimensional magic cube of a high order by constructing number sequences. Define $\Omega_n(i)$ by

$$\Omega_n(i) = \{g_1 g_2 \cdots g_i | g_k \in X, k = 1, 2, \cdots\}.$$

Because there may be two different arrangements in $\Omega_n(i)$ which correspond to the same transformation, we can get the following theorems:

***Lemma* 4.3.1** $\bigcup_{k=0}^{i} \Omega_n(k) = \bigcup_{k=0}^{i} G_n(k)$.

***Lemma* 4.3.2** $|\Omega_n(k)| \geq |G_n(k)|$.

***Lemma* 4.3.3** $\sum_{k=0}^{i} |\Omega_n(k)| \geq \sum_{k=0}^{i} |G_n(k)|$.

***Lemma* 4.3.4** Let $d$ denote the Cayley Graph diameter of $G_n$, if $\sum_{k=0}^{i} |\Omega_n(k)| < |G_n|$, then $d>i$.

***Proof.*** Use reduction to absurdity. Assume that $d \leq i$. For that $|G_n| = \sum_{k=0}^{d} |G_n(k)|$, we have

$$\sum_{k=0}^{d} |G_n(k)| > \sum_{k=0}^{i} |\Omega_n(k)| \geq \sum_{k=0}^{i} |G_n(k)|.$$

That is to say,

$$\sum_{k=0}^{d} |G_n(k)| > \sum_{k=0}^{i} |G_n(k)|,$$

which contradicts $d \leq i$. Thus, $d>i$. Proved.

*Lemma* 4.3.4 provides a method to estimate the lower bounds. It is easy to find, after eliminating duplicates, *Lemma* 4.3.4 still holds true. we can obtain more accurate lower bounds.

The accuracy of the lower bounds is decided by how many duplicates we can eliminate. If we eliminate all the duplicates, we can get the exact values of the Cayley Graph diameter.

### 4.3.2 Eliminate duplicates

In this section the superscripts of manipulations should be understood in modulo $n$. $V_i$ and $H_i$ denote $V_i^j, H_i^j (1 \le j < n)$ respectively. See the following facts:

1. $V_i^\alpha V_i^\beta = V_i^{\alpha+\beta}, H_i^\alpha H_i^\beta = H_i^{\alpha+\beta}$,
2. $V_i^\alpha V_j^\beta = V_j^\beta V_i^\alpha, H_i^\alpha H_j^\beta = H_j^\beta H_i^\alpha$,

In fact, they are two sets of the generator relations. With them we can eliminate part of the duplicates in $\Omega_n(k)$. In this section we use arrangements of manipulations to indicate a certain transformation. In particular, we use an empty arrangement "$\Phi$" to indicate the identity transformation.

***Definition* 4.3.5** A **positive sequence manipulation arrangement** is an arrangement in which any contiguous $V_iV_j$ or $H_iH_j$ satisfies $i<j$.

***Lemma* 4.3.6** Any manipulation arrangement is equivalent to a positive sequence manipulation arrangement, of which the length is not more than the original manipulation arrangement.

***Proof.*** According to fact 1, any contiguous $V_iV_i$ or $H_iH_i$ can be merged into one. According to fact 2, any contiguous $V_iV_j$ or $H_iH_j$ which satisfies $i>j$ can be swapped into $V_jV_i$ or $H_jH_i$. An equivalent positive sequence manipulation arrangement can eventually be obtained, and the length of the arrangement is not more than the original manipulation arrangement.

From *Lemma* 4.3.6, we can eliminate duplicates by selecting the positive sequence manipulation arrangements. Let $\Delta$ be the set of all positive sequence manipulation arrangements. Redefine $\Omega_n(k)$ by

$$\Omega_n(i) = \{g_1 g_2 \cdots g_i \in \Delta \,|\, g_k \in X, k=1, 2, \cdots\}.$$

In particular, when $i = 0$, set $\Omega_n(0) = \{\Phi\}$.

In order to compute the number of positive sequence manipulation arrangements, we divide $\Omega_n(i)$ into several sets on the basis of the last manipulation.

***Definition* 4.3.7** Let $v_n^k(i)$ denote the set of the positive sequence manipulation arrangements of which the last manipulation is $V_k (1 \le k \le n)$; let $h_n^k(i)$ denote the set of the positive sequence manipulation arrangements of which the last manipulation is $H_k (1 \le k \le n)$.

***Lemma* 4.3.8** $|v_n^k(i)| = |h_n^k(i)|$.

***Proof.*** By replacing all the $V$ in $v_n^k(i)$ with $H$ and replacing all the $H$ in $h_n^k(i)$ with $V$, we can construct a bijection from $v_n^k(i)$ to $h_n^k(i)$. Therefore, $|v_n^k(i)| = |h_n^k(i)|$.

***Lemma*** **4.3.9**

$$\forall i\geq 2,\ \left|v_n^k(i)\right| = (n\text{-}1)\left(\sum_{j=1}^{k\text{-}1}\left|v_n^j(i\text{-}1)\right| + \sum_{j=1}^{n}\left|h_n^j(i\text{-}1)\right|\right).$$

***Proof.*** For $\forall i\geq 2$, each positive sequence manipulation arrangement $M$ of $v_n^k(i)$, we consider its last manipulation and the former $i$-1 step manipulations. From the *Definition* 4.3.7, the last step of the positive sequence manipulation arrangement in $v_n^k(i)$ is $V_k$ (altogether $n$-1 possible cases). Its former $i$-1 manipulations constitute a positive sequence manipulation arrangement whose length is $i$-1. If its last manipulation is $V_j(k\leq j<n)$, the last two manipulations of $M$ are $V_jV_k$ which satisfies $k\leq j$. This contradicts the definition of the positive sequence manipulation arrangement. So the last manipulation must be $V_j(1\leq j<k)$ or $H_j(1\leq j\leq n)$ and the number of the positive sequence manipulation arrangements, which satisfy this condition and of which the length is $i$-1, are $\sum_{j=1}^{k\text{-}1}\left|v_n^j(i\text{-}1)\right| + \sum_{j=1}^{n}\left|h_n^j(i\text{-}1)\right|$. Besides, there are altogether $n$-1 possible cases of $V_k$, so we have $\left|v_n^k(i)\right|=(n\text{-}1)\left(\sum_{j=1}^{k\text{-}1}\left|v_n^j(i\text{-}1)\right| + \sum_{j=1}^{n}\left|h_n^j(i\text{-}1)\right|\right)$. The proof is completed.

***Theorem*** **4.3.10**

$$\forall i\geq 2,\ \left|v_n^k(i)\right| = (n\text{-}1)\left(\sum_{j=1}^{k\text{-}1}\left|v_n^j(i\text{-}1)\right| + \sum_{j=1}^{n}\left|v_n^j(i\text{-}1)\right|\right).$$

***Proof.*** *Theorem* 4.3.10 follows from *Lemma* 4.3.8 and *Lemma* 4.3.9.

Since that $v_n^j(i), h_n^j(i) (j=1, 2, \cdots, n)$ are the partitions of $\Omega_n(i)$, we have the following theorem.

***Theorem*** **4.3.11**

$$|\Omega_n(i)| = \sum_{j=1}^{n}\left|v_n^j(i)\right| + \sum_{j=1}^{n}\left|h_n^j(i)\right| = 2\sum_{j=1}^{n}\left|v_n^j(i)\right|.$$

Combining the initial values $\left|v_n^1(1)\right| = \left|v_n^2(1)\right| = \cdots = \left|v_n^n(1)\right| = n\text{-}1$, *Theorem* 4.3.10 and *Theorem* 4.3.11, we can obtain the value of $|\Omega_n(i)|$ by iterative computation. Thus, we can compute the lower bound of the Cayley Graph diameter of the 4-order two-dimensional magic cube group with *Lemma* 4.3.4.

Taking the 4-order two-dimensional magic cube as an example, with *Theorem* 4.3.10, we can get the following recurrence formulas:

$$\left|v_4^1(i)\right| = 3\left(\left|v_4^1(i\text{-}1)\right|+\left|v_4^2(i\text{-}1)\right|+\left|v_4^3(i\text{-}1)\right|+\left|v_4^4(i\text{-}1)\right|\right),$$
$$\left|v_4^2(i)\right| = 3\left(2\left|v_4^1(i\text{-}1)\right|+\left|v_4^2(i\text{-}1)\right|+\left|v_4^3(i\text{-}1)\right|+\left|v_4^4(i\text{-}1)\right|\right),$$
$$\left|v_4^3(i)\right| = 3\left(2\left|v_4^1(i\text{-}1)\right|+2\left|v_4^2(i\text{-}1)\right|+\left|v_4^3(i\text{-}1)\right|+\left|v_4^4(i\text{-}1)\right|\right),$$
$$\left|v_4^4(i)\right| = 3\left(2\left|v_4^1(i\text{-}1)\right|+2\left|v_4^2(i\text{-}1)\right|+2\left|v_4^3(i\text{-}1)\right|+\left|v_4^4(i\text{-}1)\right|\right),$$

together with initial value

$$\left|v_4^1(1)\right| = \left|v_4^2(1)\right| = \left|v_4^3(1)\right| = \left|v_4^4(1)\right| = 3,$$

we obtain the results as follows:

**Table 2 Data of the lower bound of the Cayley Graph diameter of 4-order two-dimensional magic cube**

| $i$ | $v_4^1(i)$ | $v_4^2(i)$ | $v_4^3(i)$ | $v_4^4(i)$ | $\|\Omega_n(i)\|$ | $\sum_{k=0}^{i}\|\Omega_n(k)\|$ |
|---|---|---|---|---|---|---|
| **0** | / | / | / | / | 1 | 1 |
| **1** | 3 | 3 | 3 | 3 | 24 | 25 |
| **2** | 36 | 45 | 54 | 63 | 396 | 421 |
| **3** | 594 | 702 | 837 | 999 | 6264 | 6685 |
| **4** | 9396 | 11178 | 13284 | 15795 | 99306 | 105991 |
| **5** | 148959 | 177147 | 210681 | 250533 | 1574640 | 1680631 |
| **6** | 2361960 | 2808837 | 3340278 | 3972321 | 24966792 | 26647423 |
| **7** | 37450188 | 44536068 | 52962579 | 62983413 | 395864496 | 422511919 |
| **8** | 593796744 | 706147308 | 839755512 | 998643249 | 6276685626 | 6699197545 |
| **9** | 9415028439 | 11196418671 | 13314860595 | 15834127131 | 99520869672 | 106220067217 |
| **10** | 149281304508 | 177526389825 | 211115645838 | 251060227623 | 1577967135588 | 1684187202805 |
| **11** | 2366950703382 | 2814794616906 | 3347373786381 | 3980720723895 | 25019679661128 | 26703866863933 |

Meanwhile, $|G_4|$=16!=20922789888000, for $\sum_{k=0}^{10}|\Omega_n(k)|<|G_4|\leq\sum_{k=0}^{11}|\Omega_n(k)|$, from *Lemma* 4.3.4, we can find that the Cayley Graph diameter of 4-order two-dimensional magic cube is at least 11.

According to the above discussion, by using *Python 2.7*(codes in *Appendix* 2), we obtain the results as follows (from 2-order to 8-order):

**Table 3 The comparison between estimated lower bounds and true values of the Cayley Graph diameters of different order two-dimensional magic cube groups**

| Order | Lower bound | True value |
|---|---|---|
| 2 | 3 | 4 |
| 3 | 6 | 8 |
| 4 | 11 | unknown |
| 5 | 18 | unknown |
| 6 | 26 | unknown |
| 7 | 36 | unknown |
| 8 | 48 | unknown |

# 5 The Generalization of the Model

## 5.1 Introduction

If we generalize the $n\times n$ two-dimensional magic cube as $m\times n$ form with the same rule, we can call it the $m\times n$ heteromorphic two-dimensional magic cube (we only study the cases when $m>n\geq 2$). Denoting the set of all the transformations of $m\times n$ heteromorphic two-dimensional magic cube by $G_{m\times n}$, we have the following theorem just like two-dimensional magic cube:

***Theorem* 5.1** Regard composite as operation, $G_{m\times n}$ can constitute a permutation group, and we call it the heteromorphic two-dimensional magic cube group.

## 5.2 The order of $G_{m\times n}$

In this section, we will prove that $|G_{m\times n}|=\begin{cases}(mn)^2!/2, & mn\equiv 1 \pmod 2),\\ (mn)^2!, & mn\equiv 0 \pmod 2).\end{cases}$

By making use of the method in Chapter 3, we have the following results:

***Lemma* 5.2.1** $\forall m\geq 3$, $m\in Z$, $\exists M_0\in G_{m\times n}$, $M_0=(1\ 2\ 3)$.

***Proof.*** The proof is as same as *Lemma* 3.1.

***Lemma* 5.2.2** $\forall m\geq 3$, $m\in Z$, $\{a, b, c\}\subseteq\{1, 2, \cdots, mn\}$, $\exists M_0\in G_{m\times n}$, $M_0=(a\ b\ c)$.

***Proof.*** The proof is as same as *Lemma* 3.2.

***Theorem*** **5.2.3** $G_{m\times n}=\begin{cases} A_{mn}, & mn\equiv 1(\text{mod } 2), \\ S_{mn}, & mn\equiv 0(\text{mod } 2). \end{cases}$

***Proof.*** For $G_{m\times n}$ is the permutation group applied to {1, 2, ···, *mn*}, so $G_{m\times n}\le S_{mn}$. From *Lemma* 3.2, any 3-cycle applied to {1, 2, ···, *mn*} belongs to $G_{m\times n}$. Given *Theorem* 2.1.12, all 3-cycles can generate $A_{mn}$. So $G_{m\times n}\ge A_{mn}$.

(1) When both *m* and *n* are odd, $G_n$ is generated by the even permutations $H_1^1, H_2^1, \cdots, H_n^1$ which the length is *m* and $V_1^1, V_2^1, \cdots, V_m^1$ which the length is *n*. Given *Theorem* 2.1.9, the product of any two even permutation is an even permutation, so in $G_{m\times n}$ does not exist odd permutation. Thus $G_{m\times n}\le A_{mn}$. Hence, $G_{m\times n}=A_{mn}$.

(2) When either of *m* and *n* is even, either of $H_1^1$ and $V_1^1$ must be an odd permutation, so $G_{m\times n}\ne A_{mn}$. Thus $G_{m\times n}>A_{mn}$, $|G_{m\times n}|>|A_{mn}|=(mn)!/2$ and $G_{m\times n}\le S_{mn}$. With Theorem 2.1.1, we can find that $|G_{m\times n}|$ is divisible by $|S_{mn}|=(mn)!$, thus we have $|G_{m\times n}|=(mn)!$. Therefore, $G_{m\times n}=S_{mn}$.

***Theorem*** **5.2.4** $|G_{m\times n}|=\begin{cases} (mn)^2!/2, & mn\equiv 1(\text{mod } 2), \\ (mn)^2!, & mn\equiv 0(\text{mod } 2). \end{cases}$

***Proof.*** From *Theorem* 2.1.5, *Theorem* 2.1.10 and *Theorem* 5.2.3, we can conclude that *Theorem* 5.2.4 holds true.

## 5.3 The Cayley Graph diameter of $G_{m\times n}$

The computing method of the diameter of Cayley Graph of $G_{m\times n}$ is as same as $G_n$, so we do not repeat here. Refining the program codes a little, we can compute the diameter of Cayley Graph of the heteromorphic two-dimensional magic cube groups (codes in *Appendix* 1). The results are shown as follows:

**Table 4 The Cayley Graph diameters of the partial heteromorphic two-dimensional magic cube groups**

| | 2 column | 3 column | 4 column | 5 column | 6 column |
|---|---|---|---|---|---|
| **2 row** | 4 | 7 | 8 | 12 | 14* |
| **3 row** | 7 | 8 | / | / | / |
| **4 row** | 8 | / | / | / | / |
| **5 row** | 12 | / | / | / | / |
| **6 row** | 14* | / | / | / | / |

*Note: we obtained the numbers by using the supercomputer in Ningbo University.

## 5.4 Estimated lower bound of the Cayley Graph diameter of $G_{m\times n}$

We define $v_{m\times n}^k(i)$ as the set of the positive sequence manipulation arrangements of the $m\times n$-order two-dimensional magic cubes of which the last manipulation is $V_k(1\le k\le n)$ and the length is $i$. We also define $h^k_{m\times n}(i)$ as the set of the positive sequence manipulation arrangements of the $m\times n$-order two-dimensional magic cubes of which the last manipulation is $H_k(1\le k\le n)$ and the length is $i$. It is different from Chapter 4 since that $\left|v^k_{m\times n}(i)\right|$ and $\left|h^k_{m\times n}(i)\right|$ are not equal. Through respectively discussing, we can obtain the results as follows:

$$\forall i\ge 2,\ \left|v^k_{m\times n}(i)\right|=(n\text{-}1)\left(\sum_{j=1}^{k\text{-}1}\left|v^j_{m\times n}(i\text{-}1)\right|+\sum_{j=1}^{n}\left|h^j_{m\times n}(i\text{-}1)\right|\right),$$

$$\forall i\ge 2,\ \left|h^k_{m\times n}(i)\right|=(m\text{-}1)\left(\sum_{j=1}^{k\text{-}1}\left|h^j_{m\times n}(i\text{-}1)\right|+\sum_{j=1}^{m}\left|v^j_{m\times n}(i\text{-}1)\right|\right),$$

$$|\Omega_{m\times n}(i)|=\sum_{j=1}^{n}\left|v^j_{m\times n}(i)\right|+\sum_{j=1}^{n}\left|h^j_{m\times n}(i)\right|,$$

$$\left|v^1_{m\times n}(1)\right|=\left|v^2_{m\times n}(1)\right|=\cdots=\left|v^m_{m\times n}(1)\right|=n\text{-}1,$$
$$\left|h^1_{m\times n}(1)\right|=\left|h^2_{m\times n}(1)\right|=\cdots=\left|h^n_{m\times n}(1)\right|=m\text{-}1.$$

From the above results, we can obtain the estimated lower bound of the Cayley Graph diameter of $G_{m\times n}$.

Refining the program code a little, we can compute the lower bounds of the Cayley Graph diameters of the heteromorphic two-dimensional magic cube groups (codes in *Appendix* 2). The result is given as follows:

**Table 5 The lower bounds of the diameters of the Cayley Graph of the partial heteromorphic two-dimensional magic cube groups**

| | **2 column** | **3 column** | **4 column** | **5 column** | **6 column** | **7 column** | **8 column** |
|---|---|---|---|---|---|---|---|
| **2 row** | 3 | 4 | 6 | 7 | 8 | 9 | 10 |
| **3 row** | 4 | 6 | 9 | 10 | 13 | 14 | 16 |
| **4 row** | 6 | 9 | 11 | 14 | 17 | 20 | 23 |
| **5 row** | 7 | 10 | 14 | 18 | 22 | 25 | 29 |
| **6 row** | 8 | 13 | 17 | 22 | 26 | 31 | 35 |
| **7 row** | 9 | 14 | 20 | 25 | 31 | 36 | 41 |
| **8 row** | 10 | 16 | 23 | 29 | 35 | 41 | 48 |

# 6 Open problems

**Problem 6.1** The exact values of the Cayley Graph diameters of the higher-order two-dimensional magic cube groups.

**Problem 6.2** The exact values of the Cayley Graph diameters of the heteromorphic two-dimensional magic cube groups. Are there certain rules in the Table 4?

**Problem 6.3** How to refine the lower bounds of the Cayley Graph diameters of the two-dimensional magic cube groups?

**Problem 6.4** How to compute the upper bounds of the Cayley Graph diameters of the two-dimensional magic cube groups?

## Acknowledgements

The authors would like to thank the instructor Yudong Wu for his kind instruction on the paper! And they are grateful to Professor Shihao Wei from Ningbo University and the strong support of Supercomputing Center of Ningbo University! The authors would also like to thank Tim Liu from University of Queensland and Xuansheng Wang from the University of Hong Kong for their valuable suggestions on English translation!

## Appendixes

**1. The computing program codes of calculating the diameters of Cayley Graph of *m*×*n*-orders two-dimensional magic cube groups**

```
def rotate(lst,n):
    l=len(lst)
    n=n%l
    return lst[-n:]+lst[:-n]
def move(pm,x,y,z):
    rs=list(pm[:])
    if y==0:
        temp=rotate(rs[n*(x-1):n*x:1],z)
        rs[n*(x-1):n*x:1]=temp
    else:
        temp=rotate(rs[x-1::n],z)
        rs[x-1::n]=temp
    return tuple(rs)
def god():
    global depth
    global line
    global astt
    global ct
```

```
    if len(line[depth])==0:
        depth+=1
        print "depth:"+str(depth)
        print "the number of all states:"+str(ct)
        if len(line[depth])==0:
            return False
        line.append({})
        astt[(depth-1)%3]={}
    j=line[depth].popitem()[0]
    y=0
    for x in range(1,m+1):
        for z in range(1,n):
            s=move(j,x,y,z)
            if ((s in astt[0]) or (s in astt[1]) or (s in astt[2]))==False:
                line[depth+1][s]=astt[(depth-1)%3][s]=0
                ct+=1
    y=1
    for x in range(1,n+1):
        for z in range(1,m):
            s=move(j,x,y,z)
            if ((s in astt[0]) or (s in astt[1]) or (s in astt[2]))==False:
                line[depth+1][s]=astt[(depth-1)%3][s]=0
                ct+=1
    return True
n=input("rows?")
m=input("columns?")
depth=0
ct=1
line=[{},{}]
astt=[{},{},{}]
line[0][tuple(range(0,n*m))]=[]
astt[1][tuple(range(0,n*m))]=[]
while god():
    pass
print "The God's number of two-dimensional magic cube "+str(n)+"*"+str(m)+" is "+str(depth-1)+"."
```

**2. The computing program codes of calculating the lower bounds of the diameters of Cayley Graph of *m*×*n*-orders two-dimensional magic cube groups**

```
def factorial(num):
    r=1
    for i in range(1,num+1):
        r*=i
    return r
```

```
m=input('rows?')
n=input('columns?')
v=[[n-1]*m]
h=[[m-1]*n]
d=1
s=m*(n-1)+n*(m-1)+1
print s
a=factorial(m*n)/(m*n%2+1)
while s<a:
    v.append([])
    h.append([])
    for i in range(0,n):
        temp=(n-1)*(sum(v[d-1][0:i])+sum(h[d-1]))
        v[d].append(temp)
        s+=temp
    for i in range(0,m):
        temp=(m-1)*(sum(h[d-1][0:i])+sum(v[d-1]))
        h[d].append(temp)
        s+=temp
    print s
    d+=1
print "The diameter of Cayley Graph of two-dimensional magic cube "+str(m)+"*"+str(n)+" is greater than or equal to "+str(d)+"."
```